\documentclass[12pt]{article}
\usepackage{amssymb}
\def \C{{\mathbb C}}

\def \Z{{\mathbb Z}}
\def \P{{\mathbb P}}
\def \Q{{\mathbb Q}}

\begin{document}
\title{On an automorphism of $Hilb^{[2]}$ of certain K3 surfaces}
\author{Ekaterina Amerik\thanks{Universit\'e Paris-Sud, Laboratoire de 
Math\'ematiques,
B\^atiment 425, 91405 Orsay, France; and MPIM Bonn, Vivatsgasse 7, 53111 Bonn,
Germany. e-mail: Ekaterina.Amerik@math.u-psud.fr}}

\date{July 16, 2009}
\maketitle

{\bf 1.Introduction}

\

In a recent note \cite{O}, K. Oguiso uses the structure of the cohomology
of compact hyperk\"ahler manifolds (\cite{V}) to describe the behaviour
of the dynamical degrees of an automorphism of such a manifold, and makes
an explicit computation in some particular cases. One of his examples is
the following: he considers a $K3$ surface $S$ admitting two embeddings
as a quartic in $\P^3$ (given by two different very ample line bundles $H_1$
and $H_2$). Each embedding induces an involution of the second punctual
Hilbert scheme (that is, the Hilbert scheme parametrizing finite subschemes
of length two) $X=Hilb^{[2]}S=S^{[2]}$, where a pair of general points $p_1, p_2$ is sent to
the complement of $\{p_1, p_2\}$ in the intersection of $S$ and the line 
$p_1p_2$; it is shown in
\cite{B2} that this involution is regular if and only if $S$ does not 
contain lines. Oguiso considers the product of the two involutions and shows
that this product is not induced from an automorphism of $S$, nor from
any automorphism of a K3 surface $S'$ such that $S'^{[2]}\cong S^{[2]}$.

On the other hand, in the recent past a few people have been studying the
question of potential density of rational points on $K3$ surfaces and
their symmetric powers (see for example \cite{BT}, \cite{HT}). Recall
that a variety $X$ over a number field is called potentially dense if rational
points on $X$ become Zariski-dense after a finite field extension. In 
\cite{BT}, it is proved that a $K3$ surface with an elliptic pencil
or with infinite group of automorphisms is potentially dense. The proof proceeds
by iterating rational curves by the automorphisms in the second case, and
by rational self-maps coming from the elliptic fibration (i.e. by fiberwise
multiplication by a suitable number) in the first case. In \cite{HT}, there
are several results on potential density on symmetric powers of $K3$ surfaces;
the point is that some of such symmetric products admit abelian fibrations
with a suitable potentially dense multisection, which again can be iterated.

One might ask whether the example of Oguiso leads to a new potential
density result. In fact taken as it is in \cite{O}, it does not: indeed,
Oguiso starts with a $K3$ surface $S$ of Picard number 2 whose intersection 
form represents neither
$0$ nor $-2$ on the Neron-Severi group, and for such $S$ the group of 
automorphisms is infinite (see \cite{PS}, Example in section 7). So the 
potential density is known 
already
for $S$, and a fortiori for its second punctual Hilbert scheme $X$. The 
purpose of this note is to remark that an obvious modification of Oguiso's
example gives $K3$ surfaces of Picard number 2 which carry $(-2)$-curves
(and therefore have finite automorphism group by \cite{PS}) and no elliptic
fibration, but
still admit two different embeddings as quartics in $\P^3$. For such surfaces,
potential density of the second punctual Hilbert scheme $X$ can indeed be 
proved 
using the product of the two 
involutions; whereas the intersection form on the Neron-Severi group of $X$
does not represent zero, and so there is no abelian fibration, in fact 
even no rational abelian fibration (see \cite{AC}, section 3), and thus the argument 
of \cite{HT} does not apply.

In fact Oguiso follows a remark from an earlier work of K. O'Grady
(\cite{OG}, subsection 4.4), where the author works in a much more general
situation and proposes the symplectic
manifolds equipped with two involutions satisfying certain properties
as plausible candidates for a proof of potential density. However, in the 
explicit
example given in \cite{OG}, which is $S^{[2]}$ for $S$ a general
two-dimensional linear section of the Segre embedding of
$\P^3\times \P^3$, $Aut(S)$ is again infinite by \cite{PS}, since the 
intersection
form on the Neron-Severi lattice of $S$ does not represent $0$ or $-2$.

\medskip

{\bf Acknowledgements:} I am grateful to F. Campana for sending me 
the paper \cite{O} and for stimulating discussions, and to
Max-Planck-Institut f\"ur Mathematik, Bonn, for its hospitality and support.

\

{\bf 2. The example}

\

We consider the binary quadratic form $b(x,y)=4x^2+14xy+4y^2$ (in Oguiso's
note, the form is $b'(x,y)=4x^2+16xy+4y^2$, not representing $-2$; we have
chosen ours so that it does, and the results below are also valid for
others $b_a(x,y)=4x^2+2axy+4y^2,\ a\geq 7$, representing $-2$: see the remark
at the end of this note). This is an
even indefinite form, so by \cite{M} there are K3 surfaces for which
$b$ is the intersection form on the Neron-Severi lattice. Let $S$ be such
a $K3$ surface and let $X$ be the second punctual Hilbert scheme of $S$,
so that $NS(X)\cong NS(S)\oplus \Z E$, where $E$ is one half of the class 
of the 
exceptional
divisor of the projection $X\to S^{(2)}$, where $S^{(2)}$ denotes the
symmetric square of $S$. On $NS(X)$, we have the Beauville-Bogomolov quadratic 
form $q$, defined up to a constant. It is well-known 
(\cite{B1}) that 
the direct sum decomposition is orthogonal with respect to $q$ and that for
a suitable choice of the constant, $q$ restricts as the intersection form
to $NS(S)$ and $q(E)=-2$.

\medskip

{\bf Proposition 1} {\it 1) $S$ does not carry an elliptic pencil, but
has $(-2)$-curves; in particular $Aut(S)$ is finite.

\smallskip

2) $X$ is not rationally fibered in abelian surfaces (nor in other varieties
of non-maximal Kodaira dimension).}

\medskip

{\it Proof:} The first part is immediate ($(3,-1)$ being an example of
a $(-2)$-class), except for the finiteness of $Aut(S)$ which is treated 
in \cite{PS}, section 7. For the second part, one remarks that if there is such
a fibration, then the form $q$ represents zero on $NS(X)$ (\cite{AC}). 
This is the same as to say that the equation $B(x,y)=2m^2$ has integer
solutions. But then $33y^2+8m^2$ would be a square, which is impossible
(as, for example, counting modulo 3 shows).  

\medskip

{\bf Proposition 2:} {\it On $S$, there are two classes of ample
line bundles with self-intersection 4. Those classes are very ample
and the surface $S$ does not contain lines in the correspondent projective 
embeddings}.

\medskip

{\it Proof:} Let $h_1$, $h_2$ be the classes of the line bundles corresponding
to the vectors $(1,0)$ and $(0,1)$ of the base in which we have written the intersection form.
It is immediate to check that the intersection of a nodal class (i.e. a class 
with self-intersection $-2$)
with $h_1$ or $h_2$ cannot be zero. Therefore,
using Picard-Lefschetz reflections (that is, reflections associated to the 
$(-2)$-curves) if necessary, we may assume that $h_1$
is ample (the ample cone is a chamber of the positive cone with
respect to Picard-Lefschetz reflections). To show that $h_2$ is also ample, 
it is enough to verify that the
intersection of $h_1$ with every nodal class has the same sign as the intersection
of $h_2$ with that class.

The nodal classes are $(x,y)$ satisfying $2x^2+7xy+2y^2+1=0$, so one has
$x=\frac{1}{4}(-7y\pm t)$ where $t^2=33y^2-8,\ t>0$. Now the intersection
of $(\frac{1}{4}(-7y + t), y)$  with
$H_1=(1,0)$ is equal to $t$, so we must verify that $\frac{7}{4}(-7y+t)+4y>0$,
or $7t-33y>0$. But it is immediate from $t^2=33y^2-8,\ t>0$ that
$t\geq 5y$ (and the equality holds only for $t=5,\ y=1$). So the ampleness is
proved.

The very ampleness is a consequence of the results of \cite{SD} (subsection
2.7 for the absence of base components, then Theorem 5.2 for very ampleness;
here we use the fact that $S$ does not carry an elliptic pencil).
The non-existence of lines follows by the same calculation as the ampleness.

\medskip

{\bf Corollary 3:} {\it The second punctual Hilbert scheme $X=S^{[2]}$ has two
regular involutions $\iota_1$, $\iota_2$, corresponding to the two
embeddings of $S$ by $h_1$ and $h_2$.}

\medskip

(See \cite{B2}.)

\medskip

{\bf Proposition 4:} {\it There exist $K3$ surfaces defined over a number
field whose intersection form on the Neron-Severi lattice is as above.}

\medskip

{\it Proof:} Such a $K3$ surface over $\C$ is a general member of a 
component of the Noether-Lefschetz locus of the family of quartics in $\P^3$.
Those components (consisting of quartics containing a curve of genus 3 and
degree 7) are algebraic subvarieties defined over a number field. Now the only
problem is that it could, apriori, happen that every quartic from this 
locus which
is defined over a number field has higher Picard number; but this is
ruled out by \cite{MPV}, which shows that in any family (defined over a 
number field) of smooth projective
varieties, there are members over a number field which have the same
Neron-Severi group as the general member.

\

\

{\bf 3. Potential density}
 
\

We now show that for $S$ defined over a number field, the second punctual
Hilbert scheme $X$ is potentially dense. Observe that to each embedding
of $S$ as a quartic in $\P^3$, one can associate a covering of $X$ by
a family of surfaces birational to abelian ones: in the notations of \cite{HT},
those are the surfaces $C \ast C$ where the curve $C$ runs
through the family of hyperplane sections of $S$ with one double point.
It is enough to show that the iterations of one such surface, defined over
a number field, by 
$\iota_2\iota_1$ are Zariski-dense in $X$.

Let $H_1$, $H_2$ be the elements of $NS(X)\cong NS(S)\oplus \Z E$ corresponding
to $h_1$, $h_2 \in NS(S)$ (geometrically, a divisor from  the linear system
$|H_i|$ parametrizes
subschemes whose support meets a fixed divisor from $|h_i|$).

Recall from \cite{O} that $\iota_k^*H_k=3H_k-4E$ and $\iota_k^*E=2H_k-3E$,
where $k=1,2$. Moreover the same computation as in \cite{O} gives:
$$\iota_1^*H_2=7H_1-7E-H_2,\ \iota_2^*H_1=7H_2-7E+H_1.$$ Therefore, in the
basis $\{H_1, E, H_2\}$ of $NS(X)$, the involutions $\iota_1^*$, $\iota_2^*$ are
given by the matrices 

$$M_1=\left(
\begin{array}{ccc}
3&2&7\\
-4&3&-7\\
0&0&-1
\end{array}
\right),
M_2=\left(
\begin{array}{ccc}
-1&0&0\\
-7&3&-4\\
7&2&3
\end{array}
\right).$$

The product $(\iota_2\iota_1)^*$ is thus represented by the matrix
$$M_1M_2=\left(
\begin{array}{ccc}
32&8&13\\
-24&-5&-9\\
-7&-2&-3
\end{array}
\right)$$
on $NS(X)$, and is the identity on its orthogonal complement in the second 
cohomologies of $X$ (indeed, because $h^{2,0}(X)=1$, this complement is an 
irreducible Hodge substructure, whereas $(\iota_2\iota_1)^*$ has to fix the 
holomorphic symplectic form).

\medskip

{\bf Lemma 5:} {\it No effective divisor on $X$ is invariant under 
$\iota_2\iota_1$.}

\medskip

{\it Proof:} The only divisor classes which are invariant under 
$(\iota_2\iota_1)^*$ are multiples of $L=2H_1-11E+2H_2$. These are
not effective since, for instance, the class $A=H_1-E$ is ample
(this is the inverse image of the Pl\"ucker hyperplane section by the natural 
finite morphism $X\to G(1,3)$, corresponding to the embedding of $S$ in 
$\P^3$ by $h_1$ and sending a subscheme $Z\subset S$ to the only
line in $\P^3$ containing $S$),
but its Beauville-Bogomolov intersection with $L$ is zero.   

\medskip

Let now $C_1$ be a hyperplane section of $S$ with one double point in the 
projective embedding given by $h_1$, and let $\Delta_1$ be the class
of the surface $C_1\ast C_1$ in the cohomology of $X$.

\medskip  

{\bf Proposition 6:} {\it The surface $C_1\ast C_1$ is not periodic by
$\iota_2\iota_1$.}

\medskip

{\it Proof:} It suffices to prove that $\Delta_1$ is
not periodic. Since $H^4(X, \Q)\cong S^2H^2(X, \Q)$ (from the topology
of a symmetric square) and we know the
eigenvalues of $(\iota_2\iota_1)^*$ on $H^2(X, \Q)$, one sees that
periodicity actually means invariance, so we only need to show that
$(\iota_2\iota_1)^*\Delta_1\neq \Delta_1,$ or, since $\Delta_1$ is invariant
by $\iota_1$ and $\iota_1, \iota_2$ are involutions, that 
$\iota_2^*\Delta_1\neq \Delta_1$. This in turn shall follow once we
compute that
$$\Delta_1\cdot E^2\neq \Delta_1\cdot \iota_2^*E^2=\Delta \cdot (2H_2-3E)^2.$$ 

Let $T_p\subset X$ be the surface parametrizing the length-two subschemes
of $S$ containing a given point $p$, and let $\Sigma$ be the class of
$T_p$ (obviously not depending on $p$). Since the class $H_1$ on $X$
is the class of a divisor parametrizing subschemes having some support
on the corresponding divisor on $S$, we have
$$\Delta_1 = H_1^2-q(H_1)\Sigma=H_1^2-4\Sigma.$$
Furthermore, $T_p$ is identified with the blow-up of $S$ in $p$ and
$E$ restricts to $T_p$ as the exceptional divisor, thus $\Sigma \cdot E^2=-1$.
It is also clear from the geometry that 
$$\Sigma\cdot H_1^2=\Sigma\cdot H_2^2=q(H_i)=4.$$
Moreover, for any two divisor classes $\alpha, \beta$ on $X$, one has
$$\alpha^2\cdot \beta^2=q(\alpha)q(\beta)+2q(\alpha,\beta)^2$$
(where, by abuse of notation, we denote by the same letter $q$ the 
quadratic form and the associated bilinear form), and
$$E\cdot \alpha \cdot \beta^2=0$$
(\cite{B1}).
Thus 
$$\Delta\cdot E^2=(H_1^2-4\Sigma)E^2=-8+4=-4,$$
and 
$$\Delta \cdot (2H_2-3E)^2=(H_1^2-4\Sigma)(4H_2^2-12H_2E+9E^2)>0.$$

\medskip

To sum up, we have the following

\medskip

{\bf Theorem 7:} {\it If the $K3$ surface $S$ as above is defined over
a number field, rational points are potentially dense on $X=S^{[2]}$.}

\medskip

{\it Proof:} Let $C=C_1$ be a curve as above, defined over a number field.
Since $p_g(C)=2$, $C\ast C$ is birational to an abelian surface
and hence has potentially dense rational points. Thus it suffices to show
that the union of surfaces $(\iota_2\iota_1)^k(C\ast C),\  k\in \Z$ is Zariski-dense
in $X$. By the preceding proposition, there is an infinite number of such 
surfaces, so if their union is not Zariski-dense in $X$, its Zariski closure
is a divisor. But such a divisor would be invariant by $\iota_2\iota_1$,
and this contradicts Lemma 5. 

\medskip

{\bf Remark:} In the beginning, we could have taken the binary quadratic
form $b_a(x,y)=4x^2+2axy+4y^2$ with an arbitrary $a>4$. And indeed, as soon as
this form represents $-2$, we have exactly the same results as above, up to
one exceptional case where $a=5$. In this case, $v=(1,-1)$ is a $(-2)$-class,  
and the basis vectors have intersection of different signs with $v$. Since
either $v$ or $-v$ is effective, those basis vectors cannot both represent
ample (or anti-ample) classes, even up to Picard-Lefschetz reflections. 
On the other hand, for $a=5$ (and only for $a=5$) the form $b_a$  represents 
zero as well, so that the
corresponding $K3$ surface is elliptic and therefore potentially dense.

For a general $a\geq 7$, the numbers are as follows:

$$\iota_1^*H_2=aH_1-aE-H_2;$$

$(\iota_2\iota_1)^*$ on $NS(X)$ has $1$ as an eigenvalue of multiplicity
one, the correspondent eigenvector is $2H_1-(a+4)E+2H_2$ and it is not
effective by the same reason as before. The other two eigenvalues are not
roots of unity. Thus $\Delta_1$ is invariant if periodic, and its 
non-invariance is checked in the same way as above.

\end{document}